\documentclass[a4paper,10pt]{article}

\usepackage{amssymb,amsmath}

\usepackage{natbib}  
\usepackage[OT1]{fontenc}  
\usepackage[utf8]{inputenc}   
\usepackage[english]{babel}  
\usepackage{latexsym,amssymb,amsmath,epsfig,psfrag}

\def\RR{\mathbb{R}} 
\def\EE{\mathbb{E}} 
\def\NN{\mathbb{N}}

\def\II{\mbox{ 1\hskip -.29em I}}
\newcommand{\sfrac}[2]{\kern.1em
        \raise.5ex\hbox{$#1$}\kern-.1em
        /\kern-.15em\lower.25ex\hbox{$#2$}}

\def\bdes{\begin{description}}
\def\edes{\end{description}}

\newtheorem{defi}{Definition}[section]

\newtheorem{coro}[defi]{Corollary}
\newtheorem{theo}[defi]{Theorem}
\newtheorem{exem}{Example}
\newtheorem{rem}{Remark}

\newenvironment{dem}{\vskip 2mm\noindent {\it Proof} :}
                    {\hfill $\square$ \vskip 2mm \noindent} 
 
\newcommand{\eps}{\varepsilon}
\def\bex{\begin{exem} \em }
\def\eex{\end{exem} }
\def\brem{\begin{rem} \em}
\def\erem{\end{rem} }

\begin{document}

\title{Threshold phenomena on product spaces:\\ BKKKL revisited (once more)} 

\maketitle

\begin{center}
{\sc Rapha\"el Rossignol}\footnote{Universit\'e de
  Neuch\^atel, Institut de Math\'ematiques, 11 rue Emile  Argand, Case postale
  158,  2009 Neuch\^atel, Switzerland, raphael.rossignol@unine.ch} \footnote{Rapha\"{e}l Rossignol was supported by the Swiss National Science
  Foundation grants 200021-1036251/1 and 200020-112316/1.} 
\end{center}

 \begin{abstract} 
We revisit the work of \cite{BKKKL} -- referred to as ``BKKKL'' in the title -- about influences
on Boolean functions in order to give a precise statement of threshold phenomenon on the product
space $\{1,\ldots ,r\}^\NN$, generalizing one of the main results of \cite{Talagrand94a}.\\


 \end{abstract}

\section{Introduction}
The theory of threshold phenomena can be traced back to 
\cite{Russo}, who described it as an ``approximate (Kolmogorov's) zero-one
law'' (see also \cite{Margulis} and \cite{Talagrand94a}). These phenomena occur on $\{0,1\}^n$
equipped with the probability measure $\mu_p$ which is the product of $n$
Bernoulli measures with the same parameter $p\in [0,1]$. We say that an event
$A\subset\{0,1\}^n$ is increasing if the indicator function of $A$ is
coordinate-wise nondecreasing. When the influence of each coordinate on an
increasing event $A$
is small (see the definition of $\gamma$ hereafter), and when the parameter $p$ goes from 0 to one, the probability that
$A$ occurs, $\mu_p(A)$, grows from near zero to near one on a short interval of values of $p$: this is the
\emph{threshold phenomenon}. The smaller the maximal influence
of a coordinate on $A$ is, the smaller is the bound obtained on the length of
the interval of values
of $p$. More precisely, for any $j$ in $\{1,\ldots, n\}$, define $A_j$ to be
the set of configurations in $\{0,1\}^n$ which are in $A$ and such that $j$
is \emph{pivotal} for $A$ in the following sense:
$$A_j=\{x\in\{0,1\}^n\mbox{ s.t. }x\in A,\mbox{ and }T_j(x)\not\in A\}\;,$$
where $T_j(x)$ is the configuration in $\{0,1\}^n$ obtained from $x$ by
``flipping'' coordinate $j$ to $1-x_j$. It is shown in \cite{Talagrand94a},
Corollary 1.3, that if you denote by $\gamma$ the maximum over $p$ and $j$ of the probabilities
$\mu_p(A_j)$, then, for every $p_1<p_2$,
\begin{equation}
\label{eq:Tal}
\mu_{p_1}(A)(1-\mu_{p_2}(A))\leq \gamma^{K(p_2-p_1)}\;,
\end{equation}
where $K$ is a universal constant. This result was also derived independently by \cite{FriedgutKalai}. A much simpler proof, giving the best constants up to now, was obtained later by \cite{FalikSamorodnitsky06}, and their result will be one of the main tool that we shall use in this paper. See also \cite{Rossignol06} for a more
complete description of threshold phenomena. 

This kind of phenomenon is interesting in itself, but has also
been proved useful as a theoretical tool, notably in percolation (see 
\cite{BollobasRiordan06short,BollobasRiordan06sharp,BollobasRiordan06Voronoi,Berg07}).
It seems to be partly folklore that this phenomenon occurs on other product
spaces than $\{0,1\}^n$. Notably, Theorem 5 in \cite{BollobasRiordan06Voronoi} gives a threshold result for symmetric functions on $\{1,2,3\}^n$ with an extremely short proof, mainly pointing to \cite{FriedgutKalai}. A strongly related result is \cite{BKKKL}, where it is proved
that for any subset $A$ of a product probability
space of dimension $n$, there is one coordinate that has influence of order at least
$\log n/n$ on $A$. Although the result in \cite{BKKKL} is stated in terms of influences
and not in terms of threshold phenomena, the proof can be rephrased and
slightly adapted to show
that threshold phenomena occur on various product spaces. 

Being asked by Rob van den Berg for a reference on 
generalizations of (\ref{eq:Tal}) to $\{1,\ldots ,r\}^\NN$, we could not find a truly satisfying one. The work of 
\cite{ParoissinYcart03} is close in spirit to what we were looking for, but is
stated only for symmetric sets in finite dimension. Also, Theorem 3.4 in \cite{FriedgutKalai} is even closer to what we need but is not quite adapted to $\{1,\ldots ,r\}^\NN$ since the quantity $\gamma=\max_j\mu_p(A_j)$ is replaced by the maximum of all influences, which is worse than the equivalent of $\gamma$ in $\{1,\ldots ,r\}^\NN$. The purpose of the present
note is to provide an explicit statement of the threshold phenomenon on
$\{1,\ldots ,r\}^\NN$, with a rigorous, detailed proof. We insist strongly on the fact that the  spirit
of what is written in this note can be seen as already present in \cite{BKKKL},  \cite{FriedgutKalai} and \cite{Talagrand94a}. 

Our goal will be accomplished in two steps. The first one, presented in
section \ref{sec:BKKKL}, is a general
functional inequality on the countable product
$[0,1]^\NN$ equipped with its Lebesgue measure. Then, in section \ref{sec:thresholdfinite}, we present
 the translation of this result into a threshold phenomenon on $\{1,\ldots
 ,r\}^\NN$. This is the main result of this note, stated in Corollary \ref{coro:seuilfinite}.

\section{A functional inequality on $[0,1]^\NN$, following \cite{BKKKL}}
\label{sec:BKKKL}
In \cite{Talagrand94a}, inequality (\ref{eq:Tal}) is derived from a functional inequality on $(\{0,1\}^n,\mu_p)$ (Theorem 1.5 in \cite{Talagrand94a}). Falik and Samorodnitsky's main result is also a functional inequality on $(\{0,1\}^n,\mu_p)$, with a slightly different flavour but the same spirit: it improves upon the classical Poincar\'{e} inequality essentially when the discrete partial derivatives of the function at hand have low $L^1$-norm with respect to their $L^2$-norm. Such inequalities have been extended to some continuous settings in \cite{BenaimRossignolarxiv06b}, where they were called ``modified Poincar\'e inequalities''. The discrete partial derivative is then replaced by a semi-group which is required to satisfy a certain hypercontractivity property.

In this section, we take a different road to generalize the modified Poincar\'e inequality of Talagrand (Theorem 1.5 in \cite{Talagrand94a}). This is done by combining the approach of \cite{BKKKL} and  \cite{FalikSamorodnitsky06}. This is also very close in spirit to what is done in \cite{Friedgut04}.  We will obtain  a functional inequality on $[0,1]^\NN$ equipped with the Lebesgue measure, which can be seen as a modified Poincar\'e inequality. All measures considered in this section are Lebesgue measures on Lebesgue
measurable sets.

First, we need some notations. Let $(x_{i,j})_{\substack{i\geq 1\\ j\geq 0}}$ be independent symmetric Bernoulli random
variables. For each $j$, the random variable $\sum_{i\geq
  1}\frac{x_{i,j}}{2^i}$ is uniformly distributed on $[0,1]$, whereas
$\sum_{i=1}^m\frac{x_{i,j}}{2^i}$ is uniformly distributed on
$\{\frac{k}{2^m};\;k=0,\ldots, 2^m-1\}$. For positive integers $m$ and $n$,
define a random variable $X^{m,n}$ with values in $[0,1]^\NN$ as follows:
$$(X^{m,n})_j=\left\lbrace\begin{array}{ll}\sum_{i=1}^m\frac{x_{i,j}}{2^i}&\mbox{if
    }j\leq n\\ & \\
\sum_{i\geq 1}\frac{x_{i,j}}{2^i}&\mbox{if }j> n
\end{array}
\right.$$
For any real function $f$ on $[0,1]^\NN$, we define the following random variables:
$$\Delta_{i,j}^{m,n}f= f(X^{m,n})-\EE_{x_{i,j}}\lbrack f(X^{m,n})\rbrack\;,$$
where $\EE_{x_{i,j}}$ denotes the expectation with respect to $x_{i,j}$ only.
Define $\lambda$ to be the Lebesgue measure on $[0,1]$, and if $f$
belongs to $L^2([0,1]^\NN,\lambda^{\otimes \NN})$, denote by $Var_{\lambda}(f)$
the variance of $f$ with respect to $\lambda^{\otimes \NN}$.

Finally, define, for any positive integer $n$ and any real numbers $y_1,\ldots,y_n$:
$$f_n(y_0,\ldots,y_n)=\int f(y_0,\ldots,y_n,y_{n+1},\ldots)\;\bigotimes_{k\geq
  n+1}d\lambda(y_{k})\;.$$
We shall use the following hypothesis on $f$:
\begin{equation}
\label{eq:hypriemann}
\mbox{ For every integer }n, f_n\mbox{ is Riemann-integrable.}
\end{equation}
The following result can be seen as a generalization of Theorem 1.5 in \cite{Talagrand94a}.
\begin{theo}
\label{theo:BKKKLbis}
Let $f$ be a real measurable function on $[0,1]^\NN$. Define, for $p\geq 0$:
$$N_p(f)=\limsup_{n\infty}\limsup_{m\infty}\sum_{j=0}^n\sum_{i=1}^m\EE(|\Delta_{i,j}^{m,n}f|^p)^{\frac{2}{p}}\;.$$
Suppose that $f$
belongs to $L^2([0,1]^\NN)$, and satisfies hypothesis (\ref{eq:hypriemann}). Then,
$$N_2(f)\geq \frac{1}{2}Var_{\lambda}(f)\log\frac{Var_{\lambda}(f)}{N_1(f)}\;.$$\end{theo}
\begin{dem}
Denote by $Y^{m,n}$ the first $n$ coordinates of $X^{m,n}$. Theorem 2.2 in \cite{FalikSamorodnitsky06} implies that:
$$\sum_{j=0}^n\sum_{i=1}^m\EE(\Delta_{i,j}^{m,n}f^2)\geq
\frac{1}{2}Var(f_n(Y^{m,n}))\log\frac{Var(f_n(Y^{m,n}))}{\sum_{j=0}^n\sum_{i=1}^m\EE(|\Delta_{i,j}^{m,n}f|)^2}\;.$$
Notice that $\EE(f_n(Y^{m,n}))$ is a Riemann-sum of $f_n$ over
$[0,1]^n$. Since $f_n$ is Riemann-integrable, $Var(f_n(Y^{m,n}))$ converges to
$Var(f_n(U_0,\ldots,U_n))$ when $m$ goes to infinity, where $U_0,
\ldots,U_n$ are independent random variables with uniform distribution on
$[0,1]$. Then, this is for instance a consequence of Doob's convergence
theorems for martingales bounded in $L^2$, $f_n(U_1,\ldots,U_n)$ converges in
$L^2$ to $f$ as $n$ tends to infinity. Thus,
$$\lim_{n\infty}\lim_{m\infty}Var(f_n(Y^{m,n}))=Var_{\lambda}(f)\;.$$
The theorem follows.
\end{dem}
If a function $f$
is coordinate-wise nondecreasing, we shall say it is \emph{increasing}. Now, we can get a simplified version of Theorem \ref{theo:BKKKLbis} for increasing functions. To this end, let us define the random variable $X^{\infty}$ with values in $[0,1]^\NN$ as follows:
$$\forall j\geq 0,\;(X^{\infty})_j= \sum_{i\geq 1}\frac{x_{i,j}}{2^i}\;,$$
and let:
$$\Delta_{i,j}^{\infty}f= f(X^{\infty})-\EE_{x_{i,j}}\lbrack f(X^{\infty})\rbrack\;.$$
\begin{coro}
\label{coro:croissantes}
Let $f$ be a real measurable function on $[0,1]^\NN$, increasing for the coordinate-wise
partial order. Define, for $p\geq 0$:
$$M_p(f)=\sum_{j=0}^\infty\sum_{i=1}^\infty\EE(|\Delta_{i,j}^{\infty}f|^p)^{\frac{2}{p}}\;.$$Then,
$$M_2(f)\geq \frac{1}{2}Var_{\lambda}(f)\log\frac{Var_{\lambda}(f)}{M_1(f)}\;.$$\end{coro}
\begin{dem}
We only need to show that $f$ satisfies the hypotheses of Theorem
\ref{theo:BKKKLbis}, and that $N_p(f)\leq M_p(f)$, at least when $p$ equals 1
and 2. Since $f$ is coordinate-wise increasing, so is $f_n$ for every $n$, and
thus hypothesis (\ref{eq:hypriemann}) is satisfied. The function $f$ is
trivially in $L^2([0,1]^\NN,\lambda^{\otimes \NN})$ since it is a real
measurable increasing function on $[0,1]^\NN$, and therefore is bounded. We shall use the following notation: for
$\eps\in\{0,1\}$, $f(X^{m,n}|x_{i,j}=\eps)$ denotes the value of $f$ at
$X^{m,n}$ where the value of $x_{i,j}$ is forced to be $\eps$, and for $t\in
[0,1]$, $f(X^{m,n}|y_j=t)$ denotes the value of $f$ at
$X^{m,n}$ where the value of $(X^{m,n})_j$ is replaced by $t$. We also use the
notation $\EE_{(x_{i',j})_{i'<i}}(g(X^{m,n}))$ to denote the expectation with
respect to the random variables $(x_{i',j})_{i'<i}$. For $j\leq
n$, and $p\geq 1$,
\begin{eqnarray}
\nonumber&&\EE_{(x_{i',j})_{i'\leq i}}(|\Delta_{i,j}^{m,n}f|^p)\\
\label{eq:deltanabla}&=& \frac{1}{2^p}\EE_{(x_{i',j})_{i'<i}}(|f(X^{m,n}|x_{i,j}=1)-f(X^{m,n}|x_{i,j}=0)|^p)\;,\\
\nonumber &=
&\frac{1}{2^p}\sum_{\eps\in\{0,1\}^{i-1}}\frac{1}{2^{i-1}}\left|f\left(X^{m,n}\left|y_j=\sum_{i'=1}^{i-1}\frac{\eps_{i'}}{2^{i'}}+\frac{1}{2^i}+\sum_{i'>i}\frac{x_{i',j}}{2^{i'}}\right.\right)\right.\\
\nonumber &&-\left.f\left(X^{m,n}\left|y_j=\sum_{i'=1}^{i-1}\frac{\eps_{i'}}{2^{i'}}+\sum_{i'>i}\frac{x_{i',j}}{2^{i'}}\right.\right)\right|^p\;,\\
\nonumber&=
&\frac{1}{2^{p+i-1}}\sum_{k=0}^{2^{i-1}-1}\left|f\left(X^{m,n}\left|y_j=\frac{k}{2^{i-1}}+\frac{1}{2^i}+\sum_{i'>i}\frac{x_{i',j}}{2^{i'}}\right.\right)\right.\\
\nonumber &&-\left.f\left(X^{m,n}\left|y_j=\frac{k}{2^{i-1}}+\sum_{i'>i}\frac{x_{i',j}}{2^{i'}}\right.\right)\right|^p\;,
\end{eqnarray}
Let us define $t_k=\frac{k}{2^{i}}+\sum_{i'>i}\frac{x_{i',j}}{2^{i'}}$. Notice
that $t_k<t_{k+1}$. Then,

\begin{eqnarray*}
&&\EE_{(x_{i',j})_{i'\leq i}}(|\Delta_{i,j}^{m,n}f|^p)\\
&=&\frac{1}{2^{p+i-1}}\sum_{k=0}^{2^{i-1}-1}|f(X^{m,n}|y_j=t_{2k+1})-f(X^{m,n}|y_j=t_{2k})|^p\;,\\
&\leq &\frac{(2\|f\|_{\infty})^{p-1}}{2^{p+i-1}}\sum_{k=0}^{2^{i-1}-1}|f(X^{m,n}|y_j=t_{2k+1})-f(X^{m,n}|y_j=t_{2k})|\;,\\
&\leq &\frac{\|f\|_{\infty}^p}{2^{i-1}}\;,
\end{eqnarray*}
since $f$ is increasing. Thus, when $n$ and $j$ are fixed,
$(\EE(|\Delta_{i,j}^{m,n}f|^p\II_{i\leq m})^2)_{i\geq 1}$ is dominated by
$(2^{2-2i}\|f\|_{\infty}^{2p})_{i\geq 1}$, whose sum converges. On the other hand, since $f$ is coordinate-wise increasing, the
function $y_j\mapsto f((y_i)_{i\geq 1})$ is Riemann-integrable for any fixed
$(y_i)_{i\not = j}$ and any $j$. Thus,
$$\lim_{m\infty}\EE(|\Delta_{i,j}^{m,n}f|^p)=\EE(|\Delta_{i,j}^{\infty}f|^p)\;.$$
Therefore, by Lebesgue's dominated convergence theorem,
$$
\lim_{m\infty}\sum_{i=1}^m\EE(|\Delta_{i,j}^{m,n}f|^p)^{\frac{2}{p}}=\sum_{i=1}^\infty\EE(|\Delta_{i,j}^{\infty}f|^p)^{\frac{2}{p}}\;,$$
which implies $N_p(f)=M_p(f)$. The result follows from Theorem \ref{theo:BKKKLbis}.
\end{dem}

\section{Threshold phenomenon on $\{1,\ldots,r\}^\NN$}
\label{sec:thresholdfinite}
Let $r$ be a positive integer. Let $I=]a,b[$ be a connected open subset of
$\RR$ with $a<b$, and for every $t$ in $I$, let $\mu_t$  be a probability measure on
$\{1,\ldots,r\}$, $\nu_{t,n}$ be the product measure $\mu_t^{\otimes n}$
on $H_n=\{1,\ldots,r\}^n$ and $\nu_{t,\NN}$ be the product measure $\mu_t^{\otimes \NN}$
on $H_\NN=\{1,\ldots,r\}^\NN$. We suppose that for every $k$ in $\{1,\ldots,r\}$,
the function $t\mapsto \mu_t(\{k\})$ is differentiable on $I$, and that for
every $k$ in $\{2,\ldots,r\}$, $t\mapsto \mu_t(\{k,k+1,\ldots,r\})$ is strictly
increasing. Then, we suppose that:
$$\lim_{t\rightarrow a}\mu_t(\{1\})=1,\mbox{ and }\lim_{t\rightarrow
  b}\mu_t(\{r\})=1\;.$$
The following result is a generalization of Corollary 1.3 in \cite{Talagrand94a}.
\begin{coro}
\label{coro:seuilfinite}
Let $A$ be an increasing measurable subset of $\{1,\ldots , r\}^\NN$. Let
$t_1\leq t_2$ be two  real numbers of $I$. Define:
$$\gamma_t:=\sup_j\nu_{t,\NN}(A_j)\;,$$
$$\gamma_*=\sup_{t\in [t_1,t_2]}\left\lbrace\max
\{\gamma_t,\gamma_t\log\frac{1}{\gamma_t}\}\right\rbrace\;,$$
and 
$$S^*=\inf_{t\in
  [t_1,t_2]}\inf_{k=2,\ldots ,r}\frac{d}{dt}\mu_t(\{k,k+1,\ldots,r\})\;.$$
Then,
$$\nu_{t_1,\NN}(A)(1-\nu_{t_2,\NN}(A))\leq \gamma_*^{S^*(t_2-t_1)}\;.$$

\end{coro}
\begin{dem}
Let $f=\II_A$. Suppose first that $A$ depends only on a finite number of
coordinates. Then,
$$\frac{d}{dt}\nu_{t,\NN}(A)= \sum_{j\geq 0}\int
\sum_{k=1}^r\mu_t'(k)f(x|x_j=k)\;d\nu_{t,\NN}(x)\;,$$
where $\mu_t'(k)=\frac{d}{dt}\mu_t(\{k\})$. Define, for any $k\in\{1,\ldots,r\}$,
$$S_{t,k}:=\sum_{l=k}^r\mu_t'(l)=\frac{d}{dt}\mu_t(\{k,k+1,\ldots,r\})\;.$$
By hypothesis, $S_{t,k}\geq 0$ for any $k$ in $\{2,\ldots,r\}$. Notice also that $S_{t,1}=0$. Letting $S_{t,r+1}:=0$, we have:
\begin{eqnarray*}
\sum_{k=1}^r\mu_t'(k)f(x|x_j=k)&=&\sum_{k=1}^r(S_{t,k}-S_{t,k+1})f(x|x_j=k)\;,\\
&=&\sum_{k=2}^rS_{t,k}(f(x|x_j=k)-f(x|x_j=k-1))\;.
\end{eqnarray*}
Define:
$$S_t^*=\inf_{k=2,\ldots ,r}S_{t,k}> 0\;.$$
Since $f$ is the indicator function of an increasing
  event $A$ in $H_\NN$,
$$\sum_{k=1}^r\mu_t'(k)f(x|x_j=k)\geq S_t^* (f(x|x_j=r)-f(x|x_j=1))\;.$$
Thus,
\begin{equation}
\label{eq:russofinite}
\frac{d}{dt}\nu_{t,\NN}(A)\geq S_t^*\sum_{j\geq 0}\int
f(x|x_j=r)-f(x|x_j=1)\;d\nu_{t,\NN}(x)\;.
\end{equation}
Now, we do not suppose anymore that $A$ depends on finitely many coordinates. Define, for any real function $g$ on $[0,1]$,
$$\frac{d^+}{dt} g(t)=\liminf_{t'\downarrow t}\frac{g(t')-g(t)}{t'-t}\;.$$
It is a straightforward generalization of Russo's formula for general increasing
events (see (2.28) in \cite{Grimmett2}, and the proof p.~44) to obtain from
inequality (\ref{eq:russofinite}) that when $A$ is measurable, and $f=\II_A$,:
\begin{equation}
\label{eq:russoinfinite}
\frac{d^+}{dt}\nu_{t,\NN}(A)\geq S_t^*\sum_{j\geq 0}\int
f(x|x_j=r)-f(x|x_j=1)\;d\nu_{t,\NN}(x)\;.
\end{equation}
Define $I(f)$ the total sum of influences for the event $A$:
$$I(f)=\sum_{j\geq 0}\int
f(x|x_j=r)-f(x|x_j=1)\;d\nu_{t,\NN}(x)\;.$$
Let $(u_{j})_{j\geq 0}$ be a sequence in $[0,1]^\NN$. Define a function $F_t$
from $[0,1]^\NN$ to $\{1,\ldots,r\}^\NN$ as follows:
$$\forall j\in \NN,\;\forall i\in\{1,\ldots ,r\},\;(F_t(u))_j=i\mbox{ if }\mu_t(\{1,\ldots ,i-1\})\leq
u_j <\mu_t(\{1,\ldots ,i\})\;.$$
Of course, under $\lambda^\NN$, $F_t(u)$ has distribution $\nu_{t,\NN}$. Define
$g_t$ to be the increasing, measurable function $f\circ F_t$ on $[0,1]^\NN$. Using 
Corollary \ref{coro:croissantes},
\begin{equation}
\label{eq:applicoro}
M_2(g_t)\geq
\frac{1}{2}Var_{\lambda}(g_t)\log\frac{Var_{\lambda}(g_t)}{M_1(g_t)}\;.
\end{equation}
First, notice that:
\begin{equation}
\label{eq:Varfgt}
Var_{\lambda}(g_t)=Var(f)=\nu_{t,\NN}(A)(1-\nu_{t,\NN}(A))\;.
\end{equation}
Then, according to equation (\ref{eq:deltanabla}), and since $g_t$ is
increasing and non-negative, 
\begin{eqnarray*}
\sum_{i=1}^\infty\EE(|\Delta_{i,j}^{m,n}g_t|^2)& \leq
&\frac{1}{4}\sum_{i=1}^\infty\frac{1}{2^{i-1}}\EE(g_t(X^{m,n}|y_j=1)-g_t(X^{m,n}|y_j=0))\;,\\
&=&\frac{1}{2}\int f(x|x_j=r)-f(x|x_j=1)\;d\nu_{t,\NN}(x)\;.
\end{eqnarray*}
Thus,
\begin{equation}
\label{eq:M2gt}
M_2(g_t)\leq \frac{1}{2}I(f)\;,
\end{equation}
Similarly,
\begin{equation}
\label{eq:M1gt}
\sum_{j\geq 0}\sum_{i=1}^\infty\EE(|\Delta_{i,j}^{\infty}g_t|)\leq
I(f)\;.
\end{equation}
For any $j$ in $\NN$, define $A_j$ to be the set of configurations in $\{1,\ldots
,r\}^\NN$ which are in $A$ and such that $j$ is pivotal for $A$:
$$A_j=\{x:\;x\in A,\mbox{ and }f(x|x_j=1)=0\}\;.$$
Since  $g_t$ is
increasing, 
\begin{eqnarray*}
\EE(|\Delta_{i,j}^{\infty}g_t|)&=&\EE(g_t(X^\infty)-g_t(X^\infty|x_{i,j}=0))\;,\\
&\leq &\EE(g_t(X^\infty)-g_t(X^\infty|y_j=0))\;,\\
&=&\int f(x)-f(x|x_j=1)\;d\nu_{t,\NN}(x)\;.
\end{eqnarray*}
and thus, for any $i\geq 1$,
\begin{equation}
\label{eq:gamma}
\EE(|\Delta_{i,j}^{\infty}g_t|)\leq \nu_{t,\NN}(A_j)\;.
\end{equation}
Define $\gamma_t:=\sup_j\nu_{t,\NN}(A_j)$. From (\ref{eq:M1gt})
and (\ref{eq:gamma}), we get:
$$M_1(g_t)\leq \gamma_tI(f)\;.$$
This inequality together with (\ref{eq:applicoro}), (\ref{eq:Varfgt}) and
(\ref{eq:M2gt}) leads to:
\begin{equation}
\label{eqdisj}
I(f)\geq Var(f)\log\frac{Var(f)}{\gamma_t I(f)}\;.
\end{equation}
Therefore,

$\bullet$ either $I(f)> Var(f)\log\frac{1}{\gamma_t}$,

$\bullet$ or $I(f)\leq Var(f)\log\frac{1}{\gamma_t}$, and in
this case, plugging this inequality into the right-hand side of (\ref{eqdisj}),
$$I(f)\geq Var(f)\log\frac{1}{\gamma_t\log\frac{1}{\gamma_t}}\;.$$
In any case, defining $\gamma_t^*=\sup
\{\gamma_t,\gamma_t\log\frac{1}{\gamma_t}\}$, it follows from (\ref{eq:russoinfinite}) that:
$$\frac{d^+}{dt} \nu_{t,\NN}(A)\geq S_t^*\nu_{t,\NN}(A)(1-\nu_{t,\NN}(A))\log\frac{1}{\gamma_t^*}\;.$$
Now, let $\gamma_*=\sup_{t\in [t_1,t_2]}\gamma_t^*$ and $S^*=\inf_{t\in
  [t_1,t_2]S_t^*}$. We get:
$$\frac{d^+}{dt}
\left\lbrack\log\frac{\nu_{t,\NN}(A)}{1-\nu_{t,\NN}(A)}-tS^*\log\frac{1}{\gamma_*}\right\rbrack
\geq 0\;,$$
for any $t$ in $[t_1,t_2[$.  It follows from Proposition 2, p.~19 in
\cite{Bourbaki49IV} (it is important to notice that the proof of this
Proposition works without modification if the function $f$ equals
$g+h$ where $g$ is increasing and $h$ continuous, and if the right-derivative
is replaced by $d^+/dt$) that:
\begin{eqnarray*}
\log\frac{\nu_{t_2,\NN}(A)(1-\nu_{t_1,\NN}(A))}{\nu_{t_1,\NN}(A)(1-\nu_{t_2,\NN}(A))}&\geq &(t_2-t_1)S^*\log\frac{1}{\gamma_*}\;,\\
\frac{\nu_{t_1,\NN}(A)(1-\nu_{t_2,\NN}(A))}{\nu_{t_2,\NN}(A)(1-\nu_{t_1,\NN}(A))}&\leq&\gamma_*^{S^*(t_2-t_1)}\;,
\end{eqnarray*}
and the result follows.
\end{dem}
{\bf Remark:} If one wants a cleaner version of the upperbound of Corollary \ref{coro:seuilfinite} in terms of
$\eta_*:=\sup_{t\in [t_1,t_2]}\sup_j\nu_{t,\NN}(A_j)$, simple calculus shows
that $\gamma_*\leq \eta_*^{1-1/e}\leq \eta_*^{1/2}$, which leads to:
$$\nu_{t_1,\NN}(A)(1-\nu_{t_2,\NN}(A))\leq \eta_*^{S^*(t_2-t_1)/2}\;.$$

\vspace*{1cm}

\begin{center}
{\bf Acknowledgements}
\end{center}
I would like to thank Rob van den Berg for having impulsed this writing, and for numerous helpful comments.


\begin{thebibliography}{}

\bibitem[Benaim and Rossignol, 2006]{BenaimRossignolarxiv06b}
Benaim, M. and Rossignol, R. (2006).
\newblock Exponential concentration for {F}irst {P}assage {P}ercolation through
  modified {P}oincar\'e inequalities.
\newblock \verb+http://arxiv.org/abs/math.PR/0609730+ to appear in {A}nnales de
  l'{IHP}.

\bibitem[Bollob{\'a}s and Riordan, 2006a]{BollobasRiordan06Voronoi}
Bollob{\'a}s, B. and Riordan, O. (2006a).
\newblock The critical probability for random {V}oronoi percolation in the
  plane is 1/2.
\newblock {\em Probab. Theory Related Fields}, 136(3):417--468.

\bibitem[Bollob{\'a}s and Riordan, 2006b]{BollobasRiordan06sharp}
Bollob{\'a}s, B. and Riordan, O. (2006b).
\newblock Sharp thresholds and percolation in the plane.
\newblock {\em Random Structures Algorithms}, 29(4):524--548.

\bibitem[Bollob{\'a}s and Riordan, 2006c]{BollobasRiordan06short}
Bollob{\'a}s, B. and Riordan, O. (2006c).
\newblock A short proof of the {H}arris-{K}esten theorem.
\newblock {\em Bull. London Math. Soc.}, 38(3):470--484.

\bibitem[Bourbaki, 1949]{Bourbaki49IV}
Bourbaki, N. (1949).
\newblock {\em \'{E}l\'ements de math\'ematique. {IX}. {P}remi\`ere partie:
  {L}es structures fondamentales de l'analyse. {L}ivre {IV}: {F}onctions d'une
  variable r\'eelle (th\'eorie \'el\'ementaire). {C}hapitre {I}:
  {D}\'eriv\'ees. {C}hapitre {II}: {P}rimitives et int\'egrales. {C}hapitre
  {III}: {F}onctions \'el\'ementaires}.
\newblock Actualit\'es Sci. Ind., no. 1074. Hermann et Cie., Paris.

\bibitem[Bourgain et~al., 1992]{BKKKL}
Bourgain, J., Kahn, J., Kalai, G., Katznelson, Y., and Linial, N. (1992).
\newblock The influence of variables in product spaces.
\newblock {\em Israel J. Math.}, 77:55--64.

\bibitem[Falik and Samorodnitsky, 2006]{FalikSamorodnitsky06}
Falik, D. and Samorodnitsky, A. (2006).
\newblock Edge-isoperimetric inequalities and influences.
\newblock to appear \verb+http://arxiv.org/pdf/math.CO/0512636+.

\bibitem[Friedgut, 2004]{Friedgut04}
Friedgut, E. (2004).
\newblock Influences in product spaces: {KKL} and {BKKKL} revisited.
\newblock {\em Combin. Probab. Comput.}, 13(1):17--29.

\bibitem[Friedgut and Kalai, 1996]{FriedgutKalai}
Friedgut, E. and Kalai, G. (1996).
\newblock Every monotone graph property has a sharp threshold.
\newblock {\em Proc. Amer. Math. Soc.}, 124:2993--3002.

\bibitem[Grimmett, 1999]{Grimmett2}
Grimmett, G. (1999).
\newblock {\em Percolation}, volume 321 of {\em Grundlehren der Mathematischen
  Wissenschaften [Fundamental Principles of Mathematical Sciences]}.
\newblock Springer-Verlag, Berlin, second edition.

\bibitem[Margulis, 1974]{Margulis}
Margulis, G. (1974).
\newblock Probabilistic characteristics of graphs with large connectivity.
\newblock {\em Prob. Pedarachi Inform.}, 10(2):101--108.
\newblock In {Russian}.

\bibitem[Paroissin and Ycart, 2003]{ParoissinYcart03}
Paroissin, C. and Ycart, B. (2003).
\newblock Zero-one law for the non-availability of multistate repairable
  systems.
\newblock {\em Int. J. of Reliability, Quality and Safety Engineering},
  10(3):311--322.

\bibitem[Rossignol, 2006]{Rossignol06}
Rossignol, R. (2006).
\newblock Threshold for monotone symmetric properties through a logarithmic
  {S}obolev inequality.
\newblock {\em Ann. Probab.}, 34(5):1707--1725.

\bibitem[Russo, 1982]{Russo}
Russo, L. (1982).
\newblock An approximate zero-one law.
\newblock {\em Z. Wahrscheinlichkeitstheor. Verw. Geb.}, 61:129--139.

\bibitem[Talagrand, 1994]{Talagrand94a}
Talagrand, M. (1994).
\newblock On {R}usso's approximate zero-one law.
\newblock {\em Ann. Probab.}, 22:1576--1587.

\bibitem[van~den Berg, 2007]{Berg07}
van~den Berg, J. (2007).
\newblock Approximate zero-one laws and sharpness of the percolation transition
  in a class of models including 2d ising percolation.
\newblock http://arXiv.org/pdf/0707.2077.

\end{thebibliography}

\end{document}